\pdfoutput=1
\RequirePackage{ifpdf}
\ifpdf % We~are running pdfTeX in pdf mode
\documentclass[pdftex]{sigma}
\else
\documentclass{sigma}
\fi

\usepackage[all]{xy}

\numberwithin{equation}{section}

\newtheorem{Theorem}{Theorem}[section]
\newtheorem*{Theorem*}{Theorem}
\newtheorem{Corollary}[Theorem]{Corollary}
\newtheorem{Lemma}[Theorem]{Lemma}

\theoremstyle{definition}
\newtheorem{Definition}[Theorem]{Definition}
\newtheorem*{Notation}{Notation}
\newtheorem{Example}[Theorem]{Example}
\newtheorem{Remark}[Theorem]{Remark}

\newcommand{\Set}{\mathbf{Set}}
\newcommand{\Art}{\mathbf{Art}}

\renewcommand{\bar}{\overline}
\newcommand{\de}{\partial}

\newcommand{\Oh}{\mathcal{O}}
\newcommand{\sO}{\mathcal{O}}

\newcommand{\sC}{\mathcal{C}}
\newcommand{\sD}{\mathcal{D}}
\newcommand{\sF}{\mathcal{F}}
\newcommand{\sE}{\mathcal{E}}
\newcommand{\sG}{\mathcal{G}}
\newcommand{\sH}{\mathcal{H}}

\newcommand{\sK}{\mathcal{K}}
\newcommand{\sL}{\mathcal{L}}
\newcommand{\sM}{\mathcal{M}}
\newcommand{\sP}{\mathcal{P}}
\newcommand{\sT}{\mathcal{T}}

\newcommand{\sN}{\mathcal{N}}

\newcommand{\Mor}{\operatorname{Mor}}
\newcommand{\Def}{\operatorname{Def}}
\newcommand{\Hom}{\operatorname{Hom}}

\newcommand{\HOM}{\operatorname{\mathcal H}\!{\rm om}}

\newcommand{\DER}{\operatorname{\mathcal D}\!{\rm er}}

\newcommand{\Spec}{\operatorname{Spec}}

\newcommand{\Ext}{{\operatorname{Ext}}}

\newcommand{\coker}{\operatorname{coker}}

\newcommand{\MC}{\operatorname{MC}}

\newcommand{\Z}{\mathbb{Z}}
\renewcommand{\H}{\mathbb{H}}

\newcommand{\K}{\mathbb{K} }
\newcommand\C{\mathbb{C}}
\newcommand\Del{\operatorname{Del}}
\newcommand\Tot{\operatorname{Tot}}
\newcommand\Grpd{\mathbf{Grpd}}

\begin{document}

\allowdisplaybreaks

\newcommand{\arXivNumber}{2507.20276}

\renewcommand{\PaperNumber}{009}

\FirstPageHeading

\ShortArticleName{Joint Deformations of Manifolds, Coherent Sheaves and Sections}

\ArticleName{Joint Deformations of Manifolds,\\ Coherent Sheaves and Sections}

\Author{Donatella IACONO~$^{\rm a}$ and Marco MANETTI~$^{\rm b}$}

\AuthorNameForHeading{D.~Iacono and M.~Manetti}

\Address{$^{\rm a)}$~Dipartimento di Matematica, Universit\`a degli Studi di Bari Aldo Moro, \\
\hphantom{$^{\rm a)}$}~Via E.~Orabona~4, 70125 Bari, Italy}
\EmailD{\mail{donatella.iacono@uniba.it}}

\Address{$^{\rm b)}$~Dipartimento di Matematica~Guido Castelnuovo, \\
\hphantom{$^{\rm b)}$}~Universit\`a degli studi di Roma ``La Sapienza'', P.~le Aldo Moro~5, 00185 Roma, Italy}
\EmailD{\mail{manetti@mat.uniroma1.it}, \mail{marco.manetti@uniroma1.it}}

\ArticleDates{Received August 06, 2025, in final form January 14, 2026; Published online February 03, 2026}

\Abstract{We describe a differential graded Lie algebra controlling infinitesimal deformations of triples $(X,\mathcal{F},\sigma)$, where $\mathcal{F}$ is a coherent sheaf on a smooth variety $X$ over a~field of characteristic 0 and $\sigma\in H^0(X,\mathcal{F})$. Then, we apply this result to investigate deformations of pairs (variety, divisor).}

\Keywords{deformation of manifolds; coherent sheaves and sections; differential graded Lie algebras}

\Classification{14D15; 17B70; 13D10; 14B10}

\section{Introduction}
Let $Z$ be a possibly singular hypersurface on a smooth variety $X$.
It is known, at least to experts in deformation theory,
that deformations of the pair $(X,Z)$ are controlled by the hypercohomology of the complex of sheaves
\[ \sD \colon \ \Theta_X\xrightarrow{\gamma}\sN_{Z|X} .\]
More precisely, considering the tangent sheaf $\Theta_X$ in degree $0$ and the normal sheaf
$\sN_{Z|X}$ in degree~$1$, the space of first-order deformations is~$\H^1(\sD)$, while
obstructions are contained in~$\H^2(\sD)$. The kernel of $\gamma$ is the sheaf $\Theta_X(-\log Z)$ of vector fields tangent to $Z$, that controls the locally trivial deformations of the pair
$(X,Z)$ \cite{donarendiconti,LMDT}.

The cokernel of $\gamma$ is precisely Schlessinger's $T^1$-sheaf $\sT^1_Z$ of the variety $Z$ and this controls, locally, the deformations of $Z$; note that, since $X$ is smooth, the local deformations of $Z$ are the same, up to isomorphism, as the local deformations of the pair $(X,Z)$.

While both sheaves $\Theta_X$ and $\Theta_X(-\log Z)$ carry natural structures of sheaves of DG-Lie algebras, the same does not hold for the complex $\sD$ (see next Example~\ref{ex.nonDGLie}): this may be a problem if one wants to study the deformations of the pair $(X,Z)$ via DG-Lie algebras, or more modestly if one is interested in computing the primary obstruction and the Massey products.

If $Z$ is defined as the zero locus of a section $\sigma$ of a line bundle $\sL\in \operatorname{Pic}(X)$, then the deformations of the pair $(X,Z)$ correspond to the deformations of the triple $(X, \sL, \sigma)$, up to isomorphism.

In \cite{Welters}, the author investigated the deformation theory of the triple $(X, \sL, \sigma)$. In particular, he proved that the first-order deformations are controlled by the first hypercohomology group of the complex
\begin{equation} \label{welter complex}
\sP(X,\sL)\xrightarrow{e_\sigma}\sL,
\end{equation}
where $\sP(X,\sL)$ is the sheaf of first-order differential operators of $\sL$ and $e_\sigma$ is the evaluation at~$\sigma$. The complex \eqref{welter complex} is quasi-isomorphic to the complex $\sD$, for $Z$ equal to the divisor of~$\sigma$ (see Remark \ref{remark H1 e H2 divisore}).

The main goal of this paper is to extend this analysis. First of all, we observe that the above complex \eqref{welter complex} has a natural structure of a coherent sheaf of DG-Lie algebras and we prove that the DG-Lie algebra of derived global sections controls the deformations of the triple $(X, \sL, \sigma)$ in the usual sense of Maurer--Cartan solutions/gauge equivalence.

More generally, we provide an explicit construction of a DG-Lie algebra controlling deformations of a triple $(X,\sF,\sigma)$, where $X$ is a smooth variety, $\sF$ a coherent sheaf on $X$ and $\sigma\in H^0(X,\sF)$ a global section of $\sF$ over $X$. More precisely,
taking a finite locally free resolution~$\sE^*$ of the sheaf $\sF$ and a section $s \in H^0\bigl(X,\sE^0\bigr)$ lifting $\sigma$, we are able to prove the following theorem.

\begin{Theorem*}[Theorem \ref{main theorem of triple}]
 The infinitesimal deformations of the triples $(X,\sF,\sigma)$ are controlled by the DG-Lie algebra $L$ of the derived global sections of the complex
\begin{equation}\label{equ.complexPE}
\sP^*(X,\sE^*)\xrightarrow{e_s}\sE^*.
\end{equation}
Moreover, the homotopy type of the DG-Lie algebra $L$ is independent of the choice of $\sE^*$ and $s$.
\end{Theorem*}

Our choice of considering only one section is motivated by the simplicity of exposition; it will be clear from our proof that a similar result holds also for any finite sequence of sections~${\sigma_1,\ldots,\sigma_n}$, whenever the complex~\eqref{equ.complexPE} is replaced by
\[\sP^*(X,\sE^*)\xrightarrow{(e_{s_1},\ldots,e_{s_n})}\bigoplus_{i=1}^n\sE^*.\]

The idea behind the proof of the above theorem is the following. The choice of a~pair~$(\sE^*,s)$ and of an affine open cover of $X$ provides a representative for the DG-Lie algebra
$L$ (the Thom--Whitney totalization of the \v{C}ech cochain hypercomplex), together with a (ordinary) functor~${\sG\colon \Art_\K \to \Grpd}$, such that its connected components are the same as the
isomorphism classes of the deformations of the triple $(X,\sF, \sigma)$.
Then, using the Fiorenza--Iacono--Martinengo variation of Hinich's theorem on descent of
Deligne groupoids (see Theorems~\ref{teorema di hinich} and~\ref{teorema discesa FIM}), we are able to show that both $\sG$ and the Deligne groupoid of $L$ have the same connected components.

\begin{Notation} We denote by $\K$ a field of characteristic 0, by $\Art_{\K}$ the category of local Artin $\K$-algebras with residue field $\K$, by $\Set$ the category of sets and by $\Grpd$ the (ordinary) category of groupoids.
By variety, we mean an integral (and so irreducible and reduced), separated scheme of finite type over $\K$.
\end{Notation}

\section{A brief review of Deligne groupoids of DG-Lie algebras} \label{section review Deligne groupoid}

For the readers' convenience, we review the Deligne groupoid and the deformation functor associated to a DG-Lie algebra; we refer to \cite{LMDT} for a more detailed exposition and proofs. For simplicity of exposition, we define a \emph{formal groupoid} as a covariant functor from $\Art_{\K}$ to the (ordinary) category of groupoids.

For any DG-Lie algebra $L$, the associated Deligne groupoid
is the formal groupoid
\[\Del_L\colon\ \Art_{\K}\to \Grpd\]
 defined in the following way. Given $A\in \Art_{\K}$, with maximal ideal $\mathfrak{m}_A$, fix a representative for the tensor product $L\otimes\mathfrak{m}_A$. Then the objects of $\Del_L(A)$ are the solutions of the Maurer--Cartan equation in $L\otimes\mathfrak{m}_A$
\[ \operatorname{Ob}(\Del_L(A))=\MC_L(A)=\big\{x\in L^1\otimes\mathfrak{m}_A\mid {\rm d}x+[x,x]/2=0\big\}.\]
The exponential group of the nilpotent Lie algebra $L^0\otimes\mathfrak{m}_A$ acts on $\MC_L(A)$ by the gauge action, namely
\[ e^a\ast x=x+\sum_{n\ge 0}\frac{[a,\cdot]^{n}}{(n+1)!}([a,x]-{\rm d}a),\qquad a\in L^0\otimes\mathfrak{m}_A,\qquad x\in \MC_L(A).\]
In particular, for every $\nu\in L^{-1}\otimes \mathfrak{m}_A$ and every $x\in \MC_L(A)$ we have
$e^{{\rm d}\nu+[x,\nu]}\ast x=x$ since
\[ [{\rm d}\nu+[x,\nu],x]-{\rm d}({\rm d}\nu+[x,\nu])=0.\]
Moreover, for a fixed $x\in \MC_L(A)$, the subset
$I(x):=\big\{{\rm d}\nu+[x,\nu]\mid \nu\in L^{-1}\otimes\mathfrak{m}_A\big\}$ is a Lie subalgebra of $L^0\otimes\mathfrak{m}_A$ and then $\exp(I(x))$ is a subgroup contained in the stabiliser of $x$ under the gauge action.

Given $x,y\in \MC_L(A)$, the set of morphisms $x\to y$ is the quotient set
\[ \Mor_{\Del_L(A)}(x,y)=\big\{a\in L^0\otimes\mathfrak{m}_A\mid e^a\ast x=y\big\}/{\sim},\]
where $a\sim b$ if $e^a=e^be^{{\rm d}\nu+[x,\nu]}$ for some $\nu\in L^{-1}\otimes \mathfrak{m}_A$.
As proved in \cite[Lemma 6.5.5]{LMDT}, we have $a\sim b$ if and only if there exists
$\mu\in L^{-1}\otimes \mathfrak{m}_A$ such that
$e^a=e^{{\rm d}\mu+[y,\mu]}e^b$, and this implies that the composition of morphisms is
properly defined.
Note that when $L^{-1}=0$, the Deligne groupoid is the same as the action groupoid for the gauge action on Maurer--Cartan elements.

Every morphism of DG-Lie algebras $L\to M$ naturally induces a natural transforma\-tion ${\Del_L\to \Del_M}$. Furthermore, the basic theorem of homotopy invariance of Deligne groupoids ensures that, if $L\to M$ is a quasi-isomorphism of DG-Lie algebras, then
$\Del_{L}(A)\to \Del_M(A)$ is an equivalence of groupoids, for every $A\in \Art_{\K}$.

The deformation functor
 \[\Def_L\colon\ \Art_{\K}\to \Set\]
associated to $L$ is the functor of isomorphism classes in the Deligne groupoids, i.e., $\Def_L(A)=\pi_0(\Del_L(A))$ for every $A\in \Art_{\K}$; equivalently, $\Def_L$ is given by Maurer--Cartan elements modulo gauge action.

For simplicity of exposition, and to be consistent with the existing literature about Hinich's descent theorem and its extensions, in this paper we only consider ordinary functors from
$\Art_{\K}$ to the category of groupoids.
We point out that an alternative, and for some aspects also conceptually clearer, approach
is to consider $\Grpd$ as a $(2,1)$-full subcategory of the 2-category~$\textbf{Cat}$, and
work with pseudo-functors from $\Art_{\K}$ to $\Grpd$ (see, e.g., \cite[Definition 3.10]{Vistoli}).

\begin{Lemma}\label{lem.primolemma}\samepage
Let $f\colon L\to M$ be a morphism of DG-Lie algebras, $A\in \Art_{\K}$ and let
\[\xymatrix{\Del_L(A)\ar[r]^-{g}\ar[d]_f&G\ar[d]^{k}\\
\Del_M(A)\ar[r]^-{h}&H}\]
be a diagram in $\Grpd$.
Assume that
\begin{enumerate}\itemsep=0pt
\item[$(1)$] $h\circ f=k\circ g$;
\item[$(2)$] $f\colon H^i(L)\to H^i(M)$ is surjective for $i=-1$ and injective for
$i=0$;
\item[$(3)$] $g$ is full and essentially surjective;
\item[$(4)$] $h$ is faithful.
\end{enumerate}
Then $g$ is an equivalence of groupoids.
\end{Lemma}

\begin{proof}
It is sufficient to prove that for every Maurer--Cartan element $x\in \MC_L(A)$ the map
\[g\colon\ \Mor_{\Del_L(A)}(x,x)\to \Mor_G(g(x),g(x))\]
is injective. According to \cite[Lemma 6.6.8]{LMDT}, hypothesis (2) implies that the morphism
\[ f\colon\ \Mor_{\Del_L(A)}(x,x)\to \Mor_{\Del_M(A)}(f(x),f(x))\]
is injective, while $h$ is injective on morphisms by assumption.

Hence $h\circ f$ is injective on morphisms and the first condition ensures that $k\circ g$ is also injective on morphisms.
\end{proof}

\begin{Remark}
If we work in the (2,1)-category of groupoids, the first item in the above lemma can be weakened by requiring that the diagram is $2$-commutative, i.e., there exists a 2-morphism~$\eta$ between $h\circ f$ and $k\circ g$. Indeed, this implies the existence of an invertible map $\eta _x$ for any $ x \in \MC_L(A)$, such that for any~${a \in \Mor_{\Del_L(A)}(x,x)}$, $k\circ g(a) =\eta_x^{-1} \circ (h\circ f)(a) \circ \eta_x$. Then, we can conclude that the injectivity of $h\circ f$ on morphisms ensures that $k\circ g$ is also injective on morphisms.
\end{Remark}

\begin{Example}\label{ex.irrelevantsuHom}
For later use, we explicitly describe the Deligne groupoid for some particular examples of DG-Lie algebras.

Let $R$ be a commutative $\K$-algebra and
let
\[(E^*,\de)\colon\ \cdots\xrightarrow{\de} E^{n-1}\xrightarrow{\de}E^n\xrightarrow{\de}\cdots \xrightarrow{\de}E^0\]
be a bounded above complex of projective $R$-modules.

The first DG-Lie algebra that we consider is
 \smash{$M=\Hom^*_R(E^*,E^*)$}, with bracket $[f,g]=fg-\smash{(-1)^{\bar{f} \bar{g}}gf}$
and differential
${\rm d}=[\de, \cdot]$, i.e., \smash{${\rm d}f=\de f-(-1)^{\bar{f}}f\de$}.

For any element $ x\in M^1\otimes \mathfrak{m}_A\subseteq
\Hom^1_{R}(E^*,E^*\otimes \mathfrak{m}_A)\subseteq \Hom^1_{R\otimes A}(E^*\otimes A,E^*\otimes A)$, the Maurer--Cartan equation can be written as $[\de+x,\de+x]/2=0$ and then
$ x\in \MC_M(A)$ if and only if
$(E^*\otimes A,\de+x)$ is a complex of $R\otimes A$-modules.

Since the differential $d$ is the inner derivation $[\de,-]$, by a completely standard computation, see, e.g., either \cite[Proposition~2.5.3 and Example~6.3.2]{LMDT} or \cite[Example 1.2]{seattle},
given $x,y\in M^1\otimes\mathfrak{m}_A$ and $a\in M^0\otimes\mathfrak{m}_A$ we have
\[ e^a\ast x=y\iff e^{[a,-]}(\de+x)=\de+y\iff e^a(\de+x)e^{-a}=\de+y,\]
where the last equation is expressed in the associative product of $\Hom^*_{R\otimes A}(E^*\otimes A,E^*\otimes A)$.
Thus, the relation $e^a\ast x=y$ means that $e^a\colon (E^*\otimes A,\de+x)\to (E^*\otimes A,\de+y)$ is a morphism of complexes
of $R\otimes A$-modules.
Finally, using the power series expansions of $e^x$ and $\log(1-x)$ it is easy to see that
for $ x\in \MC_M(A)$ the morphisms of complexes
\[e^{{\rm d}\nu+[x,\nu]}\colon\ (E^*\otimes A,\de+x)\to (E^*\otimes A,\de+x),\qquad
\nu\in M^{-1}\otimes \mathfrak{m}_A,\]
are exactly those that are homotopic to the identity by a $R\otimes A$-linear
homotopy that is trivial on $E^*$. Hence, if $e^a\ast x=e^b\ast x=y$, the equivalence
$a\sim b$ means that the two morphisms $e^a,e^b\colon (E^*\otimes A,\de+x)\to (E^*\otimes A,\de+y)$ are homotopic by a $R\otimes A$-linear homotopy that is trivial modulo $\mathfrak{m}_A$, see \cite{FIM}.
\end{Example}

\begin{Example}
As a second example, fix an element $s\in E^0$ \big(since $E^1=0$ we have $\de s=0$\big); let~$\tau$ be an indeterminate of degree $-1$ and consider the short exact sequence of complexes
\[ 0\to (E^*,\de)\to (E^*\oplus R\tau,\de)\to R\tau\to 0,\qquad\text{where}\quad \de\tau=s.\]
Define
$ L\subseteq \Hom^*_R(E^*\oplus R\tau,E^*\oplus R\tau)$ as the DG-Lie subalgebra of endomorphisms $f$ such that~${f(E^*\oplus R\tau)\subseteq E^*}$. Therefore,
we have isomorphisms
\[ L^i\cong M^i\oplus E^{i-1} =\Hom^i_R(E^*,E^*)\oplus E^{i-1},
\qquad f\mapsto (f_{|E^*},f(\tau)).\]
Via the above isomorphism $ L\cong M\oplus E^*[-1]$, the bracket in $ L$ becomes
\[ [(f,v),(g,w)]=\bigl([f,g],f(w)-(-1)^{\bar{f} \bar{g}}g(v)\bigr),\]
and the differential
\[ {\rm d}(f,v)=[(\de,s),(f,v)]=\bigl([\de,f],\de v-(-1)^{\bar{f}}f(s)\bigr).\]
Then, we have the following short exact sequence of DG-Lie algebras
\[ 0\to E^*[-1]\xrightarrow{\iota} L\xrightarrow{\rho} M\to 0,\]
where $E^*[-1]$ is equipped with the trivial bracket, $\rho$ is the projection and
$\iota(x)=(-1)^{i}(0,x)$ for~${x\in E^{i-1}}$.

We observe the following.
\begin{itemize}\itemsep=0pt
\item[(1)] Under the above isomorphisms, a pair $ (u,t)\in L^1\otimes\mathfrak{m}_A$, with $u\in \Hom^1_R(E^*,E^*)\otimes\mathfrak{m}_A$ and~${t\in E^0\otimes\mathfrak{m}_A}$, satisfies the Maurer--Cartan equation if and only if
$(\de+u)^2=0$.

\item[(2)] For every $ (u,t)\in L^1\otimes\mathfrak{m}_A$ and every
$b\in E^{-1}\otimes \mathfrak{m}_A$, we have
\begin{equation} \label{azione gauge (0,b)}
 e^{(0,b)}\ast (u,t)=(u,t-(\de+u)b).
 \end{equation}
\end{itemize}

\begin{Lemma}\label{lem.gaugeaffine}
In the above notation,
let $ (u,t)\in L^1\otimes\mathfrak{m}_A$ be a solution of the Maurer--Cartan equation,
$ (f,a)\in L^0\otimes\mathfrak{m}_A$ and denote
$\bigl(u',t'\bigr)=e^{(f,a)}\ast (u,t)$. Then
\[ \de+u'=e^f(\de+u)e^{-f},\qquad s+t'-e^{f}(s+t)\in \bigl(\de+u'\bigr)\bigl(E^{-1}\otimes\mathfrak{m}_A\bigr).
\]
In particular, $s+t'=e^{f}(s+t)$ in the cohomology group $H^0(E^*\otimes A,\de+u')$ and there exists~${b\in E^{-1}\otimes \mathfrak{m}_A}$ such that
\begin{equation}\label{equ.gauge2}
e^{(0,b)}\ast\bigl(u',t'\bigr)=e^{(0,b)}e^{(f,a)}\ast (u,t)=\bigl(u',e^{f}(s+t)-s\bigr).
\end{equation}
\end{Lemma}

\begin{proof} It is convenient to prove the above formulas in the DG-associative algebra
$\Hom^*_R(E^*\oplus R\tau,E^*\oplus R\tau)\otimes A$. Then,
the elements $(u,t)$, $\bigl(u',t'\bigr)$ and $(f,a)$ correspond, respectively, to the $A$-linear extension of the endomorphisms:
\begin{enumerate}\itemsep=0pt
\item[(1)] $\alpha\colon E^*\oplus R\tau\to E^*\otimes\mathfrak{m}_A$,
$\alpha_{|E^*}=u$, $\alpha(\tau)=t$;
\item[(2)] $\beta\colon E^*\oplus R\tau\to E^*\otimes\mathfrak{m}_A$, $\beta_{|E^*}=u'$, $\beta(\tau)=t'$;
\item[(3)] $\gamma\colon E^*\oplus R\tau\to E^*\otimes\mathfrak{m}_A$, $\gamma_{|E^*}=f$, $\gamma(\tau)=a$.
\end{enumerate}
By the relation $\de+\beta=e^\gamma(\de+\alpha)e^{-\gamma}$, we get
\[ s+t'=(\de+\beta)(\tau)=e^\gamma(\de+\alpha)e^{-\gamma}(\tau)=
e^\gamma(\de+\alpha)\left(\tau-\frac{e^{-\gamma}-I}{\gamma}(a)\right).\]
Denoting by $c=\frac{e^{-\gamma}-I}{\gamma}(a)\in E^{-1}\otimes\mathfrak{m}_A$, we conclude that
\[\begin{split} s+t'&=e^\gamma(\de+\alpha)(\tau-c)=e^\gamma(s+t-(\de+u)(c))=
e^f(s+t)-
e^{f}(\de+u)(c)\\
&=e^{f}(s+t)-\bigl(\de+u'\bigr)e^{f}(c),\end{split}\]
and therefore
\[s+t'-e^{f}(s+t)\in \bigl(\de+u'\bigr)\bigl(E^{-1}\otimes\mathfrak{m}_A\bigr),\qquad
e^{f}(s+t)-s=t'+\bigl(\de+u'\bigr)e^{f}(c).\]
According to \eqref{azione gauge (0,b)}, the element $b=-e^{f}(c)\in E^{-1}\otimes\mathfrak{m}_A$ satisfies~\eqref{equ.gauge2}.
\end{proof}
\end{Example}

\section{Descent of Deligne groupoids }

For later application, we review the extension given in \cite{FIM} (see Theorem \ref{teorema discesa FIM}) of Hinich's theorem on descent of Deligne groupoids \cite{Hinich} (see Theorem \ref{teorema di hinich}), which is particularly useful in the deformation theory of coherent sheaves that are not locally free.

Let
 \[
G\colon\ \xymatrix{ {{G}_0}
\ar@<2pt>[r]\ar@<-2pt>[r] & { {G}_1}
 \ar@<4pt>[r] \ar[r] \ar@<-4pt>[r] & { {G}_2}
\ar@<6pt>[r] \ar@<2pt>[r] \ar@<-2pt>[r] \ar@<-6pt>[r]&
\cdots}
\]
be a semicosimplicial groupoid with face operators $\delta_i\colon G_n\to G_{n+1}$ and cosimplicial identities~${\delta_{l}
\delta_{k}=\delta_{k+1} \delta_{l}}$ for $l\leq k$; in particular
$\delta_0\delta_0=\delta_1\delta_0$, $ \delta_0\delta_1=\delta_2\delta_0$, $
\delta_1\delta_1=\delta_2\delta_1$.

The associated total, or descent groupoid,
$\Tot(G)$ is defined in the following way:
\begin{enumerate}\itemsep=0pt
\item[(1)] the objects of $\Tot(G)$ are the pairs $(l,m)$,
with $l\in G_0$ and $m\colon \delta_{0}l\to \delta_{1}l$ a morphism in~$G_1$ such that the cocycle diagram
\[
\begin{matrix}\xymatrix{&\delta_0\delta_0l\ar@{=}[dl]\ar[dr]^{\delta_0m}&\\
\delta_1\delta_0l\ar[d]^{\delta_1m} && \delta_0\delta_1l\ar@{=}[d]\\
\delta_1\delta_1l\ar@{=}[dr] && \delta_2\delta_0l\ar[dl]_{\delta_2m}\\
&\delta_2\delta_1l&}\end{matrix}
\]
is commutative in $G_2$;

\item[(2)] the morphisms between $(l_0,m_0)$ and $(l_1,m_1)$ are the morphisms $a$ in $G_0$
between $l_0$ and $l_1$ such that the diagram
\[
\xymatrix{
\delta_{0}l_0\ar[r]^{m_0}\ar[d]_{\delta_{0}a}&\delta_{1}l_0\ar[d]^{{\delta_{1}a}}\\
\delta_{0}l_1\ar[r]^{m_1}&\delta_{1}l_1
}
\]
is commutative in $G_1$.
\end{enumerate}

The above construction is functorial; every morphism
$f\colon G\to\! H$ of semicosimplicial groupoids induces a morphism of groupoids
$f\colon \Tot(G)\to \Tot(H)$. Moreover,
if for every $n=0,1,2$ the component $f_n\colon G_n\to H_n$ is an equivalence of groupoids, then also
$f\colon \Tot(G)\to \Tot(H)$ is an equivalence of groupoids.

Given a semicosimplicial DG-Lie algebra
\[L\colon\
\xymatrix{ {L_0}
\ar@<2pt>[r]\ar@<-2pt>[r] & {L_1}
 \ar@<4pt>[r] \ar[r] \ar@<-4pt>[r] & { L_2}
\ar@<6pt>[r] \ar@<2pt>[r] \ar@<-2pt>[r] \ar@<-6pt>[r]& \cdots }\]
with face morphisms $\delta_i\colon L_{n-1}\to L_n$, $0\le i\le n$, we have the associated semicosimplicial Deligne groupoids
\[\Del_L\colon\
\xymatrix{ \Del_{L_0}
\ar@<2pt>[r]\ar@<-2pt>[r] & \Del_{L_1}
 \ar@<4pt>[r] \ar[r] \ar@<-4pt>[r] & \Del_{L_2}
\ar@<6pt>[r] \ar@<2pt>[r] \ar@<-2pt>[r] \ar@<-6pt>[r]& \cdots }\]
 that is a semicosimplicial object in the category of functor from $\Art_{\K}$ to $\Grpd$ and we can consider its descent groupoid $\Tot(\Del_L)\colon \Art_{\K}\to \Grpd$.

Let us denote by $\Tot(L)$ the Thom--Whitney totalization of the semicosimplicial DG-Lie algebra $L$. We refer to \cite{FIM,LMDT} for its explicit construction; here, we just recall some basic properties that will be used in this paper.
\begin{enumerate}\itemsep=0pt
\item[(1)] $\Tot(L)$ is a DG-Lie algebra and $H^*(\Tot(L))\cong H^*\bigl(\Tot^{\Pi}C(L)\bigr)$, where
\[ C(L)\colon\ \cdots\xrightarrow{\delta} L_{n-1}\xrightarrow{\delta=\sum_{i=1}^n(-1)^{i-1}\delta_i}L_n\xrightarrow{\delta}\cdots\]
is the cochain double complex associated to $L$ \cite[Section 1.2.6]{weibel}.
\item[(2)] $\Tot$ is an exact functor from the category of semicosimplicial DG-Lie algebras to the category of DG-Lie algebras.
\end{enumerate}

Here, we recall two useful results on descent of Deligne groupoids.
	\begin{Theorem}[\cite{Hinich}]\label{teorema di hinich}
			Let $L$ be a semicosimplicial DG-Lie algebra such that every $L_i$ is concentrated in nonnegative degrees. Then there exists an equivalence of groupoids $\Del_{\Tot(L)}\cong \Tot(\Del_L)$.
		\end{Theorem}

 This means that there exists a natural transformation
of functors $\Phi\colon \Del_{\Tot(L)}\to \Tot(\Del_L)$ that is an equivalence of groupoids for every Artin ring.
 This was also generalised to an equivalence to the formal groupoid associated to the $L_\infty$-structure on the homotopy limit, in the case of semicosimplicial Lie algebra \cite{FMM}.
If some $L_i$ has nontrivial elements of negative degree the above result is generally false. However, we have the following result.

\begin{Theorem}[\cite{FIM}]\label{thm.fim}\label{teorema discesa FIM}
In the above situation, if $H^p(L_i)=0$ for every $i$ and every $p<0$, then there exists an isomorphism of functors of Artin rings
$\Phi\colon \Def_{\Tot(L)}\to \pi_0(\Tot(\Del_L))$.
\end{Theorem}

\section{Coherent DG-Lie algebroids on smooth varieties}\label{sezione costruzione cocone C}

Let $X$ be a smooth variety over $\K$ (i.e., a regular, integral, separated scheme of finite type) and denote by $\Theta_X=\DER_{\K}(\Oh_X,\Oh_X)$ its tangent sheaf.

\begin{Definition}\label{def.liealgebroid}
	A \emph{coherent DG-Lie algebroid} over $X$ is a pair $(\sL^*,\alpha)$, where
	\begin{itemize}\itemsep=0pt \item $\sL^*$ is a sheaf of DG-Lie algebras over $\K$ on $X$ such that the underlying
	sheaf of DG-vector spaces
	\smash{$\cdots \sL^i\xrightarrow{\delta}\sL^{i+1}\cdots$} is a bounded complex of coherent $\sO_X$-modules (hence $\delta$ is $\sO_X$ linear, while the bracket is assumed only $\K$-bilinear);
		\item $\alpha\colon\sL^*\to\Theta_X$ is a morphism of complexes of $\sO_X$-modules, called the \emph{anchor map}, commuting with the brackets;
		\item finally, we require the Leibniz rule to hold $[ l , fm ] = \alpha(l)(f)m + f[l,m]$, $\forall l,m\in\sL^*$, $f\in\sO_X$.
	\end{itemize}
\end{Definition}

A morphism (resp.\ quasi-isomorphism) $\sL^*\to \sM^*$ of coherent DG-Lie algebroids over $X$
is a morphism (resp.\ quasi-isomorphism) of complexes of $\Oh_X$-modules commuting with brackets and anchor maps.

An example is given by the sheaf $\sP^*(X,\sE^*)$ of first-order differential operators with scalar symbol of a bounded complex $(\sE^*,\de)$ of locally free sheaves on $X$. For applications to deformation theory it is convenient, following \cite{DMcoppie}, to define
$\sP^*(X,\sE^*)$ in terms of derivations of pairs.

Consider the subcomplex of $\sO_X$-modules
\[ \begin{split}
\sD^*(X,\sE^*)&\subseteq \Theta_X \times \HOM^*_{\K}(\sE^*,\sE^*),\\[4pt]
\sD^*(X,\sE^*)&=\{(h,u)\mid u(ax)-au(x)=h(a)x, \text{ for every } a\in \Oh_X, \, x\in \sE^*\}.
\end{split}\]

Then, we define $\sP^*(X,\sE^*)\subseteq\HOM^*_{\K}(\sE^*,\sE^*)$
as the image of $\sD^*(X,\sE^*)$ under the projection on the second factor; since $\sE^*$ is torsion-free, we have an isomorphism of sheaves of DG-Lie algebras~${\sP^*(X,\sE^*)\xrightarrow{\simeq}\sD^*(X,\sE^*)}$, $ u\mapsto (\alpha(u),u)$,
where $\alpha(u)$ is the (scalar) symbol of $u$. If~${u\in \sP^n(X,\sE^*)}$ and $n\not=0$, then
\smash{$u\in \HOM^n_{\Oh_X}(\sE^*,\sE^*)$} and $\alpha(u)=0$, while every
$u\in \sP^0(X,\sE^*)$ is a sequence of differential operators of first order~${u_i\colon \sE^i\to \sE^i}$ with the same scalar symbol.

Since the complex $\sE^*$ is assumed bounded, it is not difficult to prove that $\sP^*(X,\sE^*)$ is also a coherent DG-Lie algebroid with surjective anchor map $\alpha$ \cite[Proposition 5.1]{DMcoppie}. Therefore, we~have a short exact sequence
\begin{equation}\label{eq atyah sequence for complex}
 0 \to \HOM^*_{\Oh_X}(\sE^*,\sE^*) \to \sP^*(X,\sE^*) \xrightarrow{ \alpha } \Theta_{X}\to 0.
\end{equation}

Next, we describe another example of DG-Lie algebroid, extending the previous construction, starting from a coherent sheaf and a section.

\begin{Definition}\label{def.resolution}
Let $\sF$ be a coherent sheaf on a smooth variety $X$ and $\sigma\in H^0(X,\sF)$ a global section. A resolution of $(\sF,\sigma)$ over $X$ is a pair
$(\sE^*,s)$, where
\begin{itemize}\itemsep=0pt
\item $\sE^*$ is a finite locally free resolution of $\sF$, i.e., we have an exact sequence
$0\to \sE^{-m} \smash{\xrightarrow{\de}} \cdots \smash{\xrightarrow{\de}}\allowbreak \sE^0 \to \sF \to 0$ with every $\sE^i$ locally free;
\item $s\in H^0\bigl(X,\sE^0\bigr)$ is a section lifting $\sigma$.
\end{itemize}
\end{Definition}

Since $X$ is a smooth variety, a resolution as above always exists by Kleiman's theorem, see, e.g., \cite[Exercise III.6.8]{Hartshorne} or \cite[Theorem 3.3]{Borelli}, and, for any of them, evaluation at $s$ gives a~morphism of complexes of $\sO_X$-modules
$ e_s\colon \sP^*(X,\sE^*)\to \sE^*$, $ e_s(u)=u(s)$.

Mimicking the construction of the DG-Lie algebra $M$ of Section~\ref{section review Deligne groupoid}, the mapping cocone of the morphism $e_s$, denoted by $\sC^*(X,\sE^*,s)$, carries a structure of DG-Lie algebroid. More explicitly, $\sC^n(X,\sE^*,s)=\sP^n(X,\sE^*)\oplus\sE^{n-1}$ for every
$n$, then we extend the bracket of $\sP^*(X,\sE^*)$ to $\sC^*(X,\sE^*,s)$ by setting
\[ \begin{cases}
[u,x]=u(x)&\text{for }u\in \sP^*(X,\sE^*), \ x\in \sE^*,\\
[x,y]=0&\text{for } x,y\in \sE^*,\end{cases}\]
and the differential is as follows:
\[ {\rm d}\colon\ \sC^*(X,\sE^*,s)\to \sC^*(X,\sE^*,s),\qquad
{\rm d}(u,x)=([\de,u],\de x+[s,u])=[(\de,s),(u,x)].\]
It is straightforward to show that $\sC^*(X,\sE^*,s)$ is a sheaf of DG-Lie algebras. Moreover,
$\sC^*(X,\allowbreak\sE^*,s)$ inherits the anchor map $\alpha$ of \eqref{eq atyah sequence for complex} so that it has a structure of DG-Lie algebroid.

Note that, for every $r\in H^0\bigl(X,\sE^{-1}\bigr)$, there exists an isomorphism of DG-Lie algebroids
\[% \label{eq. C^* +dr}
\sC^*(X,\sE^*,s)\xrightarrow{\simeq} \sC^*(X,\sE^*,s+\de r),\qquad (u,x)\mapsto
(u,x-[r,u]).
\]

The projection $\sC^*(X,\sE^*,s)\to \sP^*(X,\sE^*)$ is a morphism of DG-Lie algebroids and there is a~commutative diagram of sheaves of DG-Lie algebras, with exact rows and columns{\samepage
\begin{gather}\label{diagramma successioni esatte con DGLA}
\begin{matrix}\xymatrix{ & 0\ar[d]&0\ar[d] & & \\
 & \sE^*[-1] \ar@{=}[r]\ar[d] &\sE^*[-1] \ar[d] & & \\
 0 \ar[r] & \sK \ar[r] \ar[d] &\sC^*(X,\sE^*,s)\ar[r] \ar[d] & \Theta_X \ar[r]\ar@{=}[d] & 0 \\
 0 \ar[r] & \HOM^*_{\Oh_X}(\sE^*,\sE^*) \ar[d] \ar[r] & \sP^*(X,\sE^*) \ar[r]^-{\alpha} \ar[d] & \Theta_X \ar[r] &0 \\
 & 0 & 0, & & & \\
 }\end{matrix}
 \end{gather}
where $\sK$ is the mapping cocone of the evaluation map $e_s\colon \HOM^*_{\Oh_X}(\sE^*,\sE^*) \to\sE^*$.}

We point out that if $(\sE^*,s)$ is a resolution of $(\sF,\sigma)$ and
\begin{equation}\label{equ.morfismorisoluzioni}
\begin{matrix}\xymatrix{ \sE^* \ar^{\varphi}[r] \ar[dr] & \sH ^* \ar[d] \\
 & \sF }\end{matrix}
\end{equation}
is a morphism of finite locally free resolutions of $\sF$, then
$(\sH^*,\varphi(s))$ is also a resolution of the pair $(\sF,\sigma)$.

\begin{Lemma} In the above situation \eqref{equ.morfismorisoluzioni}, if $\varphi$ is injective and its cokernel is a complex of locally free sheaves, then the DG-Lie algebroids $\sC^*(X,\sE^*,s)$ and
$\sC^*(X,\sH^*,\varphi(s))$ are quasi-isomorphic.
\end{Lemma}

\begin{proof} The proof is essentially the same given in \cite{DMcoppie} for the analogous result for the DG-Lie algebroids $\sD^*(X,\sE^*)$ and $\sD^*(X,\sH^*)$.
Define
\[ \sP^*\bigl(X,\sE^*\xrightarrow{\varphi} \sH^*\bigr)\subseteq \sP^*(X,\sE^*)\times_{\Theta_X}
\sP^*(X,\sH^*)\]
as the subset of pairs $(u,v)$ such that $\varphi u=v\varphi$. It is clearly a DG-Lie algebroid on $X$ and,
according to \cite[Corollary 3.3]{DMcoppie}, the projections on both factors give two short exact sequences of complexes of $\Oh_X$-modules
\begin{gather*}
0\to \sP^*\bigl(X,\sE^*\xrightarrow{\varphi} \sH^*\bigr)\to \sP^*(X,\sH^*)\to
\HOM^*_{\Oh_X}(\sE^*,\coker\varphi)\to 0,\\
 0\to \HOM^*_{\Oh_X}(\coker\varphi,\sH^*)\to \sP^*\bigl(X,\sE^*\xrightarrow{\varphi} \sH^*\bigr)\to \sP^*(X,\sE^*)\to 0.
 \end{gather*}
Since $\varphi$ is an injective morphism of resolutions, its cokernel is acyclic; since $\sE^*$ is locally free and bounded above, the complex
$\HOM^*_{\Oh_X}(\sE^*,\coker\varphi)$ is acyclic. Similarly, since
$\coker\varphi$ is acyclic, locally free and bounded above, also the complex
$\HOM^*_{\Oh_X}(\coker\varphi,\sH^*)$ is acyclic.

Therefore, we have a commutative diagram
\[
\xymatrix{&\sP^*\bigl(X,\sE^*\xrightarrow{\varphi} \sH^*\bigr)\ar[ld]\ar[rd]&\\
\sP^*(X,\sE^*)\ar[d]^{e_s}&&\sP^*(X,\sH^*)\ar[d]^{e_{\varphi(s)}}\\
\sE^*\ar[rr]^\varphi&&\sH^*,}
\]
where the oblique arrows are quasi-isomorphisms of DG-Lie algebroids.
Denoting by $\sC^*\bigl(X,\sE^*\xrightarrow{\varphi} \sH^*,s\bigr)$ the mapping cocone of the evaluation map \smash{$\sP^*\bigl(X,\sE^*\xrightarrow{\varphi} \sH^*\bigr)\to \sE^*$}, defined as $(u,v)\mapsto u(s)$, we have a span of quasi-isomorphisms of DG-Lie algebroids
\[ \sC^*(X,\sE^*,s)\xleftarrow{\qquad}\sC^*\bigl(X,\sE^*\xrightarrow{\varphi} \sH^*,s\bigr)
\xrightarrow{\qquad}\sC^*(X,\sH^*,\varphi(s)).\tag*{\qed} \]\renewcommand{\qed}{}
\end{proof}

\begin{Theorem}\label{DGLA non dipende da risoluzione}
Let $(X,\sF,\sigma)$ be as above and
 $ (\sE^*_1,s_1)$ and
 $ (\sE^*_2,s_2)$ two resolutions of $(\sF,\sigma)$.
 Then, the DG-Lie algebroids $\sC^*(X,\sE_1^*,s_1)$ and
 $\sC^*(X,\sE_2^*,s_2)$ are quasi-isomorphic, i.e., they are connected by a zig-zag of quasi-isomorphisms of DG-Lie algebroids.
\end{Theorem}

\begin{proof} Again, the proof is essentially the same of
\cite[Lemma 7.7]{DMcoppie}. Denote by $f_1\colon \sE_1^*\to \sF$ and~${f_2\colon \sE_2^*\to \sF}$ the two surjective quasi-isomorphisms of complexes such that $f_i(s_i)=\sigma$.
Since~$X$ is smooth, the usual killing cycles procedure gives two bounded
 complexes $\sH ^*$ and $ \sN ^*$ of locally free sheaves and a commutative diagram
\[\xymatrix{0 \ar[r]&\sE_1^* \oplus \sE_2^* \ar[r]^-{j_1+j_2} \ar[dr]_{f_1+ f_2} & \sH ^* \ar[d]^{h} \ar[r] & \sN ^* \ar[r] &0 \\
 & & \sF &\\ }\]
such that $h\colon\sH ^* \to \sF$ is a resolution, the
upper row is a short exact sequence and there exists~${b\in H^0\bigl(X,\sH^{-1}\bigr)}$ such that
${\rm d}b=j_1(s_1)-j_2(s_2)$.
By previous lemmas, the DG-Lie algebroid~${\sC^*(X,\sE_i^*,s_i)}$ is quasi-isomorphic to
$\sC^*(X,\sH^*,j_i(s_i))$, for $i=1,2$, and the two DG-Lie algebroids
$\sC^*(X,\sH^*,j_1(s_1))$ and $\sC^*(X,\sH^*,j_2(s_2))$ are isomorphic.
\end{proof}

 \begin{Definition}\label{definizione comologia T i }
Let $(X,\sF,\sigma)$ be as above and
 $ (\sE^*,s)$ a resolution of $(\sF,\sigma)$. For any $i\in \Z$, we~define the coherent sheaves \smash{$ \sT^i_{ (X, \sF, \sigma)}$} as the cohomology sheaves of $\sC^*(X,\sE^*,s)$, i.e.,
$\smash{\sT^i_{(X,\sF,\sigma)}}:= \sH^i( X, \sC^*(X,\sE^*,s)) $.
Analogously, we can define the hyper-cohomology groups $\smash{T^i_{(X,\sF,\sigma)} } = \H^i( X,\allowbreak\sC^*(X,\sE^*,s))$.
By the previous Theorem~\ref{DGLA non dipende da risoluzione}, the sheaves \smash{$\sT^i_{(X,\sF,\sigma)}$} and the groups \smash{$T^i_{(X,\sF,\sigma)}$} are well defined, since they do not depend on the choice of the resolution. Moreover, the sheaves $\smash{\sT^i_{(X,\sF)}}=\sH^i( X, \sP^*(X,\sE^*) )$ and the groups \smash{$T^i_{(X,\sF)}=\H^i( X, \sP^*(X,\sE^*) )$} are also well defined.
 \end{Definition}

\begin{Remark}\label{remark T^i negativi svaniscono}
Since $\sE^*$ is a locally free resolution of the coherent sheaf $\sF$, we have
\[ \H^i(X,\sE^*[-1])=H^{i-1}(X,\sF),\qquad \H^i( X, \HOM^*_{\Oh_X}(\sE^*,\sE^*))=\Ext^i_X(\sF,\sF).\]
Therefore, it follows from \eqref{diagramma successioni esatte con DGLA} that the sheaves $\sT^i_{(X,\sF,\sigma)}$ and the groups $T^i_{(X,\sF,\sigma)}$ vanish for $i<0$, while for $i\ge 0$ we have the
long exact sequences
\[
\xymatrixcolsep{10pt}\xymatrix{
 0 \ar[r] & \H^0( X,\sK ) \ar[r] \ar[d] &T^0_{(X,\sF,\sigma)}\ar[r] \ar[d] & H^0(X,\Theta_X) \ar[r]\ar@{=}[d] & \H^1(X,\sK) \ar[r] \ar[d] &T^1_{(X,\sF,\sigma)} \ar[r] \ar[d] & \cdots\\
 0 \ar[r] & \Ext^0_X(\sF,\sF) \ar[r] & T^0_{(X,\sF)}\ \ar[r] & H^0(X,\Theta_X) \ar[r] & \Ext^1_X(\sF,\sF) \ar[r] & T^1_{(X,\sF)} \ar[r] & \cdots, \\
 }
\]
and
\begin{equation}\label{successione esatta Ti}
\xymatrixcolsep{12pt}\xymatrix{
 \cdots \ar[r] &T^i_{(X,\sF,\sigma)}\ar[r] & T^i_{(X,\sF)} \ar[r] & H^i(X,\sF) \ar[r] &T^{i+1}_{(X,\sF,\sigma)}\ar[r] & T^{i+1}_{(X,\sF)} \ar[r] & \cdots. \\
 }
\end{equation}
\end{Remark}

\begin{Remark}\label{C in caso locally free solo 0 e 1}
If $\sF$ is a locally free sheaf on $X$ and $\sigma\in H^0(X,\sF)$ a global section, then the pair $(\sF,\sigma)$ is itself a resolution.
Then, the evaluation at $\sigma$ gives a morphism of sheaves of $\sO_X$-modules~${e_\sigma \colon \sP(X,\sF)\to \sF}$, $ e_\sigma(u)=u(\sigma)$, that is, the differential of the DG-Lie algebroid~${\sC^*(X,\sF,\sigma)}$. Note that in this case the DG-Lie algebroid is concentrated in degree zero and one: $\sC^0(X,\sF,\sigma)=\sP(X,\sF)$ and $\sC^1(X,\sF,\sigma)=\sF$.
\end{Remark}

\section[Deformation of triples (X, sF,sigma)]{Deformation of triples $\boldsymbol{ (X, \sF,\sigma)}$}\label{section def triple}

In this section, we prove that, for any resolution $ (\sE^*,s)$ of a triple $(X,\sF,\sigma)$, the sheaf of DG-Lie algebras $\sC^*(X,\sE^*,s)$
 controls the deformations of the triple $(X,\sF,\sigma)$. This means that the DG-Lie algebra
$R\Gamma(\sC^*(X,\sE^*,s))$
 of derived global sections of
$\sC^*(X,\sE^*,s)$ controls the deformations in the usual sense of Maurer--Cartan solutions modulo gauge action.

\begin{Definition}
Let $X$ be a variety over $\K$, $\sF$ a~coherent sheaf on $X$ and $\sigma\in H^0(X,\sF)$ a~global section of $\sF$ over $X$. An infinitesimal deformation of $(X, \sF,\sigma)$ over $A\in \Art_\K$ is the data~${(\mathcal{X}_A, \sF_A, \sigma_A)}$, where
\begin{itemize}\itemsep=0pt
\item $\mathcal{X}_A$ is an infinitesimal deformation of $X$ over $A$, i.e., a pull-back diagram
\begin{center}
$\xymatrix{ X \ar[r] \ar[d] & \mathcal{X}_A \ar[d]^\pi \\
 \Spec \K \ar[r] & \Spec A,\\ }$
\end{center}
where $\pi$ is flat; since $A$ is an Artin ring, this is given by a sheaf $\Oh_{\mathcal{X}_A}$ of flat $A$-algebras over $X$ together with an isomorphism of sheaves \smash{$\Oh_{\mathcal{X}_A}\otimes_A\K\to \Oh_X$}.
\item $ \sF_A$ is a coherent sheaf of $\Oh_{\mathcal{X}_A}$-modules, flat over $A$, equipped with a morphism $\sF_A \to \sF$, inducing an isomorphism \smash{$\sF_A \otimes_{\Oh_{\mathcal{X}_A}} \Oh_X \to \sF$}.

\item $\sigma_A$ is a global section of $\sF_A$ extending $\sigma$.
\end{itemize}
\end{Definition}

\begin{Definition} \label{def isomorfismo defomrazioni triple}An isomorphism of (infinitesimal) deformations $(\mathcal{X}_A, \sF_A,\sigma_A)\to\bigl(\mathcal{X}_A', \sF_A',\sigma_A'\bigr)$ is the data of
\begin{enumerate}\itemsep=0pt
\item[(1)] an isomorphism of sheaves of $A$-algebras \smash{$\theta\colon \Oh_{X_{A}}\to \Oh_{X_{A}'}$} extending the identity on $\Oh_X$;
\item[(2)] an isomorphism of sheaves of $\Oh_{X_{A}}$-modules $\varphi\colon \sF_A\to \sF_A'$, extending the identity on $\sF$, where the $\Oh_{X_{A}}$-module structure on $\sF_A'$ is induced by $\theta$;
\item[(3)] $\varphi(\sigma_A)=\sigma_{A}'$.
\end{enumerate}
\end{Definition}

\begin{Definition} We denote by
\smash{$\Def_{(X,\sF, \sigma)} \colon \Art_\K \to \Set$}
the functor of isomorphism classes of infinitesimal deformations of the triple $(X, \sF,\sigma)$.
\end{Definition}

In \cite{DMcoppie},
we considered the functor $\Def_{(X,\sF)}$ of deformations of the pair $(X, \sF )$, defined in a~completely similar way. It is plain that there exists a natural forgetful natural transformation~${\Def_{(X,\sF, \sigma)}\to \Def_{(X,\sF)}}$.

\subsection*{The case $\boldsymbol{ X}$ affine}
 Assume $X$ smooth affine variety and
let \smash{$0\to \sE^{-m} \xrightarrow{\de} \cdots \xrightarrow{\de} \sE^0 \to \sF \to 0$} be a finite free resolution of $\sF$. Denote by
$M=\Gamma(X,\sD^*(X,\sE^*))$
the DG-Lie algebra of global sections of the sheaf~${\sD^*(X,\sE^*)\cong\sP^*(X,\sE^*)}$.

According to \cite[Proposition 7.9]{DMcoppie}, there exists an isomorphism of deformation functors ${\Def_M\cong \Def_{(X,\sF)}}$; for our application, we need to describe a suitable factorization of this isomorphism (see Lemma~\ref{lem.defo coppie in gruppoidi}).

For every
$A\in \Art_{\K}$, fix a representative for the tensor product $\sE^0\otimes A$ and
define the formal groupoid
\smash{$\sG_{(X,\sF)}\colon \Art_{\K}\to \Grpd$}
in the following way. For every
$A\in \Art_{\K}$, the objects of~$\sG_{(X,\sF)}(A)$ are the
$A$-flat quotients $\sF_A$ of $\sE^0\otimes A$ such that the projection map
$p\colon \sE^0\otimes A\to \sE^0$ factors to a morphism
$p\colon\sF_A\to\sF$ and induces an isomorphism $\sF_A\otimes_A\K=\sF$.

A morphism $\sF_A\to \sF_A'$ is a pair $(\theta,\varphi)$, where
\begin{enumerate}\itemsep=0pt
\item[(1)] $\theta$ is an isomorphism of sheaves of $A$-algebras $\theta\colon \Oh_X\otimes A\to \Oh_X\otimes A$, extending the identity on $\Oh_X$;
\item[(2)] $\varphi$ is an isomorphism of sheaves of $A$-modules $\varphi\colon \sF_A\to \sF_A'$, extending the identity on $\sF$ such that $\varphi(gx)=\theta(g)\varphi(x)$ for every
$g\in \Oh_X\otimes A$, $x\in \sF_A$.
\end{enumerate}

\begin{Lemma}\label{lem.defo coppie in gruppoidi}
In the above notation, there exist
an equivalence of formal groupoids
\[ \Psi\colon \Del_M\to \sG_{(X,\sF)}\]
and an isomorphism of functors of Artin rings
\smash{$\pi_0(\sG_{(X,\sF)})\simeq \Def_{(X,\sF)}$}.
\end{Lemma}

\begin{proof}
An object in $\Del_{M}(A)$ is an element $u\in \Hom^1_{\Oh_X}(\sE^*,\sE^*)\otimes \mathfrak{m}_A$ that satisfies the Maurer--Cartan equation, so that $(\de+u)^2=0$ and this defines the following complex of sheaves
\[ 0\to \sE^{-m}\otimes A \xrightarrow{\de+u} \cdots \xrightarrow{\de+u} \sE^0\otimes A \xrightarrow{q}\sF_A \to 0,\]
where $q$ is just the cokernel of \smash{$\sE^1\otimes A\xrightarrow{\de+u} \sE^0\otimes A$}.
By flatness and lifting of relations \cite[Proposition 3.1]{Ar}, the above complex is an exact sequence,
the coherent sheaf $\sF_A $ is flat over $A$ and
$\sF_A\otimes_A\K=\sF$.
Thus $\Psi(u):=\sF_A$ is an object of $\sG_{(X,\sF)}$.

Assume now that $e^{f}\ast u=u'$
is a gauge equivalence between two Maurer--Cartan elements, with $f\in M^0\otimes\mathfrak{m}_A$. Here $f\in H^0\bigl(X,\sP^0(X,\sE^*)\bigr)\otimes\mathfrak{m}_A$ is a first-order differential operator with scalar symbol $\alpha(f)\in H^0(X,\Theta_X)\otimes\mathfrak{m}_A$. Then
$e^{\alpha(f)}\colon \Oh_X\otimes A\to \Oh_X\otimes A$ is an automorphism of sheaves of $A$-algebras and it is proved in \cite{DMcoppie} that the isomorphism of complexes of sheaves of $A$-modules
$e^f\colon (\sE^*\otimes A,\de+u)\to (\sE^*\otimes A,\de+u')$ satisfies the condition
$e^f(gx)=e^{\alpha(f)}(g)e^f(x)$, for every $g\in \Oh_X\otimes A$.
In particular, $e^f$ induces an isomorphism in cohomology, i.e., $e^f\colon \sF_A\to \sF_A'$. If $f$ is of type $f={\rm d}\nu+[u,\nu]$ for some $\nu\in M^{-1}\otimes\mathfrak{m}_A$, then
$\alpha(f)=0$, i.e., $f$ is $\Oh_X$-linear and the isomorphism of complexes
$e^f$ induces the identity in cohomology.

In conclusion, we have defined a morphism of groupoids $\Psi\colon \Del_M(A)\to \sG_{(X,\sF)}(A)$, for any~${A\in \Art_\K}$.

Since $X$ is smooth affine, every deformation of $X$ is trivial. By flatness and lifting of relations every deformation of $\sF$ is obtained as above; this implies that the morphism of groupoids $\Psi\colon \Del_M(A)\to \sG_{(X,\sF)}(A)$ is surjective on objects and that the natural map $\pi_0(\sG_{(X,\sF)}(A))\to \Def_{(X,F)}(A)$ is surjective.
Since the morphisms in $\sG_{(X,\sF)}$ are, by definition, the same as the morphisms of deformations, the map $\pi_0(\sG_{(X,\sF)}(A))\to \Def_{(X,\sF)}(A)$ is also injective.

Every morphism in $\sG_{(X,\sF)}(A)$ lifts to an $A$-automorphism of the pair $(\Oh_X\otimes A, \sE^* \otimes A)$ lifting the identity on $(\Oh_X,\sE^*)$, and
according to \cite[Lemma 2.10 and Proposition 7.9]{DMcoppie} the group of such
$A$-automorphisms is naturally isomorphic to $\exp\bigl(H^0\bigl(X,\sP^0(X,\sE^*)\bigr)\otimes\mathfrak{m}_A\bigr)$. As in \cite{FIM} and in Example~\ref{ex.irrelevantsuHom},
two elements in $\exp\bigl(H^0\bigl(X,\sP^0(X,\sE^*)\bigr)\otimes\mathfrak{m}_A\bigr)$ give the same morphism in the Deligne groupoid if and only if the corresponding $A$-automorphisms of $(\Oh_X\otimes A, \sE^* \otimes A)$ are the same in cohomology, i.e., if they give the same isomorphism of deformations.
Then, $\Psi$ is also full and faithful.
\end{proof}

We now generalize the above construction introducing the formal groupoid
\[\sG_{(X,\sF, \sigma)}\colon\ \Art_{\K}\to \Grpd,\]
depending on the resolution $(\sE^*,s)$, and defined in the following way. For every
$A\in \Art_{\K}$, the objects of $\sG_{(X,\sF, \sigma)}(A)$ are the pairs
$(\sF_A,\sigma_A)$, where
\begin{enumerate}\itemsep=0pt
\item[(1)] $\sF_A$ is an
$A$-flat quotient of $\sE^0\otimes A$ such that the projection map
$p\colon \sE^0\otimes A\to \sE^0$ factors to a morphism
$p\colon\sF_A\to\sF$ and induces an isomorphism $\sF_A\otimes_A\K=\sF$,
\item[(2)] $\sigma_A\in H^0(X,\sF_A)$ is a section such that $p(\sigma_A)=\sigma$.
\end{enumerate}

A morphism $(\sF_A,\sigma_A)\to (\sF_A',\sigma_A')$ is a pair $(\theta,\varphi)$, where
\begin{enumerate}\itemsep=0pt
\item[(1)] $\theta$ is an isomorphism of sheaves of $A$-algebras $\theta\colon \Oh_X\otimes A\to \Oh_X\otimes A$, extending the identity on $\Oh_X$;
\item[(2)] $\varphi$ is an isomorphism of sheaves of $A$-modules $\varphi\colon \sF_A\to \sF_A'$, extending the identity on $\sF$ such that $\varphi(gx)=\theta(g)\varphi(x)$ for every
$g\in \Oh_X\otimes A$, $x\in \sF_A$, and $\varphi(\sigma_A)=\sigma_A'$.
\end{enumerate}

There exists an obvious forgetful map $\sG_{(X,\sF, \sigma)}\to \sG_{(X,\sF)}$ and it immediately follows from Lemma~\ref{lem.defo coppie in gruppoidi} that, on $X$ smooth affine, every deformation of the triple over $A\in \Art_{\K}$ is isomorphic to an element of
$\sG_{(X,\sF, \sigma)}(A)$; by definition, the morphisms in $\sG_{(X,\sF, \sigma)}(A)$ are the same as the morphisms between deformations of the triple.

\begin{Lemma}\label{lem.stost} Let $X$ be a smooth affine variety, $\sF$ a coherent sheaf on $X$ and $\sigma\in H^0(X,\sF)$. Given a resolution $(\sE^*,s)$ of the triple $(X,\sF,\sigma)$, let $L$ be the DG-Lie algebra of global sections of $\sC^*(X,\sE^*,s)$. Then, there exist a natural equivalence of formal groupoids
\[ \Phi\colon \Del_{L}\to \sG_{(X,\sF, \sigma)}\]
and an isomorphism of functors of Artin rings
\smash{$\pi_0(\sG_{(X,\sF,\sigma)})\simeq \Def_{(X,\sF,\sigma)}$}.
\end{Lemma}

\begin{proof}
Suppose $A\in \Art_{\K}$ and $(u,t)\in \MC_L(A)$; then
$u\in \Hom_{\Oh_X}^1(\sE^*,\sE^*) \otimes\mathfrak{m}_A$, $t\in H^0\bigl(X,\sE^0\bigr)\allowbreak\otimes\mathfrak{m}_A $ and
\begin{equation}\label{equ.complessodeformato}
0\to \sE^{-m}\otimes A \xrightarrow{\de+u} \cdots\xrightarrow{\de+u} \sE^{-1}\otimes A\xrightarrow{\de+u} \sE^0\otimes A
\end{equation}
is a complex. Since $\sE^*$ is exact in negative degrees, by standard results in
homological algebra (e.g., \cite[Corollary 4.1.2]{LMDT}), the complex~\eqref{equ.complessodeformato} is also exact in negative degrees.

Denoting by $\sF_A$ the cokernel of the rightmost map, we have an exact sequence
\[
0\to \sE^{-m}\otimes A \xrightarrow{\de+u} \cdots \xrightarrow{\de+u} \sE^{-1}\otimes A\xrightarrow{\de+u} \sE^0\otimes A\xrightarrow{q}\sF_A\to 0
\]
and the pair $\Phi(u,t):=(\sF_A,q(s+t))$ is an object of $\sG_{(X,\sF, \sigma)}(A)$;
by flatness and relation theorem, every object of \smash{$\sG_{(X,\sF, \sigma)}(A)$} is obtained in this way.

Assume now that $e^{(f,a)}\ast(u,t)=\bigl(u',t'\bigr)$
is a gauge equivalence between two Maurer--Cartan elements, with $(f,a)\in L^0\otimes\mathfrak{m}_A$. Here $f\in H^0(X,\sP^*(X,\sE^*)) \otimes\mathfrak{m}_A$ is a first-order differential operator with scalar symbol $\alpha(f)\in H^0(X,\Theta_X)\otimes\mathfrak{m}_A$. Then
$e^{\alpha(f)}\colon \Oh_X\otimes A\to \Oh_X\otimes A$ is an automorphism of sheaves of $A$-algebras and
we have proved in \cite{DMcoppie} that the isomorphism of graded sheaves of $A$-modules
$e^f\colon \sE^*\otimes A\to \sE^*\otimes A$ satisfies the condition
$e^f(gx)=e^{\alpha(f)}(g)e^f(x)$ for every $g\in \Oh_X\otimes A$.

We consider two exact sequences
\[ \sE^{-1}\otimes A\xrightarrow{\de+u} \sE^0\otimes A\xrightarrow{q}\sF_A\to 0,\qquad
\sE^{-1}\otimes A\xrightarrow{\de+u'} \sE^0\otimes A\xrightarrow{q'}\sF_A'\to 0.\]
Arguing as in the proof of Lemma~\ref{lem.gaugeaffine}, we have that $q'\bigl(e^f(s+t)\bigr)=q'(s+t')$ and then
 the pair \smash{$\Phi(f,a)=\bigl(e^{\alpha(f)},e^f\bigr)$} is a morphism between
$\Phi(u,t)$ and $\Phi\bigl(u',t'\bigr)$. If $(f,a)={\rm d}\nu+[u,\nu]$ for some~${\nu\in L^{-1}\otimes\mathfrak{m}_A}$,
 then $\alpha(f)=0$ and $e^f$ induce the identity in cohomology, hence $\Phi(f,a)$ is the identity on $\Phi(u,t)$.

Next, we prove that $\Phi$ is surjective on morphisms.
Given $(u,t),\bigl(u',t'\bigr)\in \MC_L(A)$ and an isomorphism
$(\theta,\varphi)\colon (\Oh_X\otimes A, \sF_A)\to (\Oh_X\otimes A, \sF_A')$ such that
$\varphi(q(s+t))=q'(s+t')$, by Lemma \ref{lem.defo coppie in gruppoidi} there exists
$f\in L^0\otimes\mathfrak{m}_A$ such that $e^f\ast u=u'$ and
$q'e^f=\varphi q$.
Denoting by~${(u',h)=e^{(f,0)}\ast (u,t)}$ we have seen that $q'(s+h)=q'\bigl(e^f(s+t)\bigr)$, hence also
$q'(s+h)=q'(s+t')$ and there exists $a\in H^0\bigl(X,\sE^{-1}\bigr)$ such that
$e^{(0,a)}\ast (u',h)=\bigl(u',t'\bigr)$ (see equation \eqref{azione gauge (0,b)}).

The last step is to prove that $\Phi$ is an equivalence of groupoids. This follows by applying
Lemmas~\ref{lem.defo coppie in gruppoidi} and~\ref{lem.primolemma} to the diagram
\[\xymatrix{\Del_{L}(A)\ar[r]^-{\Phi}\ar[d]&\sG_{(X,\sF, \sigma)}(A)\ar[d]\\
\Del_{M}(A)\ar[r]^-{\Psi}&\sG_{(X,\sF)}(A),}\]
where the first vertical arrow is induced by the natural projection $L\to M$, whose kernel is the complex of vector spaces $\Gamma(X,\sE^*[-1])$, that has cohomology concentrated in degree $+1$.
\end{proof}

It is plain that the morphism $\Phi$ of Lemma~\ref{lem.stost} commutes with the restrictions to open affine subsets.

\subsection*{The general case} Let $(X,\sF,\sigma)$ be a triple with $X$ smooth
variety and $(\sE^*,s)$ a resolution of $(\sF,\sigma)$.
For notational simplicity, we shall denote by $\sC^*$ the DG-Lie algebroid $\sC^*(X,\sE^*,s)$.

Let $\mathcal{U}=\{U_i\}_i$ be any affine open cover of~$X$, then we have a semicosimplicial DG-Lie algebra
\[ \sC^*(\mathcal{U})\colon\ \xymatrix{ {\prod_i\mathcal{C}^*(U_i)}
\ar@<2pt>[r]\ar@<-2pt>[r] & { \prod_{i,j}\mathcal{C}^*(U_{ij})}
 \ar@<4pt>[r] \ar[r] \ar@<-4pt>[r] &
 {\prod_{i,j,k}\mathcal{C}^*(U_{ijk})}
\ar@<6pt>[r] \ar@<2pt>[r] \ar@<-2pt>[r] \ar@<-6pt>[r]& \cdots},\]
where the face operators \[ \partial_{h}\colon
{\prod_{i_0,\ldots ,i_{k-1}} \sC^*(U_{i_0 \dots i_{k-1}})}\to
{\prod_{i_0,\ldots ,i_k} \sC^*(U_{i_0 \dots i_k})}
\]
are given by
\[
\partial_{h}(x)_{i_0 \ldots i_{k}}={x_{i_0 \ldots
\widehat{i_h} \ldots i_{k}}}_{|U_{i_0 \dots i_k}} \qquad \text{for $ h=0,\ldots, k$}.
\]

We denote by $\Tot(\mathcal{U}, \sC^*)$ the Thom--Whitney totalization of the semicosimplicial DG-Lie algebra $\sC^*(\mathcal{U})$.
The quasi-isomorphism class of $\Tot(\mathcal{U}, \sC^*)$ is independent of
the choice of the affine cover \cite{FIM} and only depends on the quasi-isomorphism class of $\sC^*$ as coherent DG-Lie algebroid, which does not depend on the choice of the resolution by Theorem \ref{DGLA non dipende da risoluzione}. In order to prove that
$\Tot(\mathcal{U}, \sC^*)$ controls the deformations of the triple $(X,\sF,\sigma)$, we may assume that the restriction of the resolution $\sE^*$ to every $U_i$ is a finite free resolution of $\sF_{|U_i}$.

\begin{Theorem} \label{main theorem of triple}
Let $X$ be a smooth variety over $\K$, $\sF$ a coherent sheaf on $X$, $\sigma\in H^0(X,\sF)$, $(\sE^*,s)$ a resolution of $(X,\sF,\sigma)$. Let $\mathcal{U}$
be an affine open cover of $X$ as above and $\Tot(\mathcal{U}, \sC^*)$ the Thom--Whitney totalization of the semicosimplicial DG-Lie algebra $\sC^*(\mathcal{U})$, as above.
 Then, there exists an isomorphism of deformation functors
\smash{$\Def_{\Tot(\mathcal{U}, \sC^*)} \cong \Def_{(X,\sF, \sigma)}$}.
In particular, in the notation of Definition~{\rm\ref{definizione comologia T i }}, $T^1_{(X,\sF,\sigma)} $ is the tangent space and $T^2_{(X,\sF,\sigma)} $ is an obstruction space.
\end{Theorem}

 \begin{proof}
 Consider the semicosimplicial formal groupoids
 \[ \Del_{\sC^*(\mathcal{U})}\colon\ \xymatrix{ {\prod_i\Del_{\sC^*(U_i)}}
\ar@<2pt>[r]\ar@<-2pt>[r] & { \prod_{i,j}\Del_{\sC^*(U_{ij})}}
 \ar@<4pt>[r] \ar[r] \ar@<-4pt>[r] &
 {\prod_{i,j,k}\Del_{\sC^*(U_{ijk})}}
\ar@<6pt>[r] \ar@<2pt>[r] \ar@<-2pt>[r] \ar@<-6pt>[r]& \cdots},\]
 associated to the semicosimplicial DG-Lie algebra $\sC^*\!(\mathcal{U})\hspace{-0.26pt}$ and the
 semicosimplicial formal groupoid
\[ \sG(\mathcal{U})\colon\ \xymatrix{ {\prod_i\sG_{(U_i,\sF,\sigma)}}
\ar@<2pt>[r]\ar@<-2pt>[r] & { \prod_{i,j}\sG_{(U_{ij},\sF,\sigma)}}
 \ar@<4pt>[r] \ar[r] \ar@<-4pt>[r] &
 {\prod_{i,j,k}\sG_{(U_{ij},\sF,\sigma)}}
\ar@<6pt>[r] \ar@<2pt>[r] \ar@<-2pt>[r] \ar@<-6pt>[r]& \cdots}.\]

By Lemma~\ref{lem.stost}, there exists an equivalence of semicosimplicial formal groupoids
\[ \Phi\colon \Del_{\sC^*(\mathcal{U})}\to \sG(\mathcal{U}).\]
By the usual glueing procedure of schemes and sheaves, the functor
$\Def_{(X,\sF,\sigma)}$ is isomorphic to~$\pi_0(\Tot(\sG(\mathcal{U}))$.

The equivalence $\Phi$ induces an isomorphism
$\pi_0(\Tot(\sG(\mathcal{U})) \cong \pi_0(\Tot( \Del_{\sC^*(\mathcal{U})}))$.
By Remark~\ref{remark T^i negativi svaniscono}, the involved DG-Lie algebras have no negative cohomology; then we can apply Theorem~\ref{thm.fim} to conclude that
 $\pi_0(\Tot( \Del_{\sC^*(\mathcal{U})}))\cong
 \Def_{\Tot(\mathcal{U}, \sC^*)}$.

 In conclusion, we have proved a sequence of three isomorphisms of deformation functors
 \[\Def_{(X,\sF,\sigma)}\cong \pi_0(\Tot(\sG(\mathcal{U})) )\cong \pi_0(\Tot( \Del_{\sC^*(\mathcal{U})}))\cong
 \Def_{\Tot(\mathcal{U}, \sC^*)}.\tag*{\qed}
 \]\renewcommand{\qed}{}
 \end{proof}

\begin{Remark}
The case of infinitesimal deformations of a locally free sheaf together with a~subspace of sections, but on a~fixed variety, together with application through Brill--Noether theory, was already analysed using DG-Lie algebra and formal groupoids in \cite{DE-Brill-Noether} and \cite{DE-map-locally-free}.
\end{Remark}

Finally, we are able to recover a well known fact.

\begin{Corollary}\label{cor H^1=0 allora forgetful smooth}
Let $X$ be a smooth variety over $\K$, $\sF$ a coherent sheaf on $X$ and $\sigma\in H^0(X,\sF)$. If $H^1(X, \sF)=0$, then the forgetful morphism
\smash{$
\pi\colon \Def_{(X,\sF,\sigma)} \to \Def_{(X,\sF)}
$}
is smooth.

\end{Corollary}

\begin{proof}
According to the exact sequence \eqref{successione esatta Ti}, if $H^1(X, \sF)=0$, then $\pi$ is surjective on tangent spaces and injective on obstructions. Then, applying the standard smoothness criterion, $\pi$ is smooth \cite[Theorem 3.65]{LMDT}.
\end{proof}

\section{Deformation of pairs (manifold, divisor)}

Let $X$ be a smooth variety over $\K$ and $Z\subset X$ a divisor, i.e., the schematic zero locus of a~global section $\sigma$
of a line bundle $\sL$; by standard notation we can also write $\sL=\Oh_X(Z)$
and we have~${\sN_{Z / X}=\sL_{|Z}}$.
Consider
 the exact sequences
\[
0\to \Oh_X\xrightarrow{\sigma}\sL\xrightarrow{\beta}\sN_{Z|X}=\HOM_{\Oh_X}(I_Z,\Oh_Z)\to 0,
\]
and
\smash{$
0 \to {\Theta}_X(-\log Z) \to {\Theta}_X\xrightarrow{\gamma} \sN_{Z / X} \to \mathcal{T}^1_Z \to 0$},
where $\gamma $ is the obvious restriction map from ${\Theta}_X=\DER( \sO_X,\sO_X)\to \HOM_{\Oh_X}(I_Z,\Oh_Z)$. Therefore, the kernel is the subsheaf $ {\Theta}_X(-\log Z) $ of derivation of $\sO_X$ preserving the ideal sheaf $I_Z$, while
the cokernel of $\gamma$ is the Schlessinger's $ \mathcal{T}^1_Z$ sheaf which is supported on the singular locus of $Z$ (see \cite[Section 6]{Ar} and \cite{Sernesi}).

The sheaf $\sL$ is a locally free sheaf of rank~1 and so the exact sequence in
\eqref{eq atyah sequence for complex} reduces to the Atiyah extension
$0 \to \Oh_X \to \sP(X,\sL) \to \Theta_X \to 0$.

The evaluation morphism $e_\sigma \colon \sP(X,\sL)\to \sL$ defined in
Section~\ref{sezione costruzione cocone C} fits
in the following commutative diagram:
\[ \xymatrix{0\ar[r]&\Oh_X\ar[r]\ar@{=}[d]&\sP(X,\sL)\ar[d]^{e_\sigma}\ar[r] &\Theta_X\ar[r]\ar[d]^\gamma&0\\
0\ar[r]&\Oh_X\ar[r]^\sigma&\sL\ar[r]^{\beta}&\sN_{Z|X}\ar[r]&0.}\]
In fact, fixing a local nonvanishing section $\psi$ of $\sL$, we have $\sigma=f\psi$, where
$f$ is a local generator of $I_Z$. The same choice of $\psi$ allows to write every differential operator $\eta\in \sP(X,\sL)$ as $g+\chi$, with $g\in \Oh_X$, $\chi\in \Theta_X$ and $e_{\sigma}(\eta)=(gf+\chi(f))\psi$; in particular, both
$\beta(e_{\sigma}(\eta))$ and $\gamma(\chi)$ are the morphism $I_Z\to \Oh_Z$ that sends $f$ to $\chi(f)$.

By the snake lemma, we have the following commutative diagram of exact sequence:
 \[\xymatrix{ & & & &0\ar[d] & 0\ar[d] & \\
 & & & & \ker e_\sigma \ar@{=}[r] \ar[d] & \Theta_X(-\log Z) \ar[d] & \\
 & & 0 \ar[r] & \Oh_X \ar[r] \ar@{=}[d] & \sP(X,\sL)\ar[r] \ar[d]^{e_\sigma} &
 \Theta_X \ar[d]^{\gamma}\ar[r] & 0 \\
 & & 0 \ar[r] & \Oh_X \ar[r]^\sigma & \sL\ar[r] \ar[d] & \sN_{Z|X}\ar[d] \ar[r] & 0 \\
 & & & & \coker e_\sigma \ar[d] \ar@{=}[r] & \mathcal{T}^1_Z \ar[d] & \\
 & & & & 0 &0. & & }\]
Note that $e_s$ and $\gamma$ are morphisms of $\Oh_X$-modules and the two complexes of sheaves
$ \sP(X,\sL)\xrightarrow{e_\sigma}\sL$, \smash{$ \Theta_X\xrightarrow{\gamma} \sN_{Z|X}$},
are quasi-isomorphic.

\begin{Remark}\label{remark H1 e H2 divisore}
It is well known, at least to experts, that the tangent space to the deformations of the pair $(X,Z)$ is $\H^1( \Theta_X \to \sN_{Z|X})$ and the obstructions lie in $\H^2( \Theta_X \to \sN_{Z|X})$ (see, e.g., \cite{Fri15}). According to the diagram above, the same holds when we replace the complex
$\gamma\colon \Theta_X \to \sN_{Z|X}$ by the complex $e_\sigma \colon \sP(X,\sL)\to \sL$.
Note that the mapping cocone of $e_s$ carries a natural structure of DG-Lie algebras, controlling deformations of the triple $(X,\sL=\Oh_X(Z),\sigma)$, while the mapping cocone of $\gamma$ does not have a DG-Lie algebras structure compatible with the Lie bracket on $\Theta_X$ (see next example).
\end{Remark}

\begin{Example}\label{ex.nonDGLie} We want to prove that the complex
$\sD=\big[\Theta_X\xrightarrow{\gamma}\sN_{Z|X}\big]$ does not carry any bracket making the
projection map $\sD\to \Theta_X$ a DG-Lie morphism.

We consider for simplicity the case
$X=\C$, with linear coordinate $t$, and $Z=\{t=0\}$, keeping in mind that the same
argument works for every reduced divisor in a smooth manifold.

At the level of global sections, the complex of sheaves $\Theta_X\xrightarrow{\gamma}\sN_{Z|X}$ becomes
\[ \C[t]\de_t\xrightarrow{\gamma}\C,\qquad \de_t=\frac{{\rm d}}{{\rm d}t},\qquad \gamma(p(t)\de_t)=p(0).\]

Assume, by contradiction,
that the bracket on vector fields
\smash{$ [t^n\de_t,t^m\de_t]=(m-n)t^{n+m-1}\de_t$}
extends to a DG-Lie structure on the whole complex and denote $\xi=\gamma(\de_t)$. It follows that
\begin{equation}\label{equ.brachetto}
[\xi,t^n\de_t]=\begin{cases}0&\text{for }n\not=1,\\
\xi&\text{for }n=1.\end{cases}\end{equation}
In fact, if $n\ge 2$, Leibniz rule implies
$0=\gamma([\de_t,t^n\de_t])=[\xi,t^n\de_t]+[\de_t,0]$,
while for $n=1$ we have~${\xi=\gamma[\de_t,t\de_t]=[\xi,t\de_t]+[\de_t,0]}$.
If $n=0$, let $a\in \C$ be such that $[\xi,\de_t]=a\xi$; then
we deduce
\[ [\xi,\de_t]=[\xi,[\de_t,t\de_t]]=[[\xi,\de_t],t\de_t]+[\de_t,[\xi,t\de_t]]=
a\xi-a\xi=0.\]
In conclusion, the equalities \eqref{equ.brachetto} contradict Jacobi identity, since
\[ 0=\big[\xi,\big[\de_t,t^2\de_t\big]\big]=2[\xi,t\de_t]=2\xi.\]
\end{Example}

Finally, consider the DG-Lie algebroid $\sC^*(X,\sL,\sigma)$ associated to $e_\sigma\colon \sP(X,\sL)\to\sL$, that has~$e_\sigma$ as differential and it is concentrated in degree zero and one. As in the previous section, for any affine open cover $\mathcal{U}$ of $X$, we can consider the semicosimplicial DG-Lie algebra $\sC^*(\mathcal{U})=\sC^*(X,\sL,\sigma)(\mathcal{U})$.

\begin{Corollary} \label{corollario defo (X,Z)}
Let $X$ be a smooth variety over $\K$ and $Z\subset X$ a closed subscheme defined as the zero locus of the section $\sigma$
of the line bundle $\sL=\Oh_X(Z)$. Let $\mathcal{U}$
be an affine open cover of $X$ and $\Tot(\mathcal{U}, \sC^*)$ the Thom--Whitney totalization of the semicosimplicial DG-Lie algebra~${\sC^*(X,\sL,\sigma)(\mathcal{U})}$.
 Then, there exist two isomorphisms of deformation functors
 \[\Def_{\Tot(\mathcal{U}, \sC^*)} \cong \Def_{(X,\sL, \sigma)} \cong \Def_{(X,Z)}.\]
\end{Corollary}

\begin{proof}
By general facts in deformation theory \cite[Section 3.4.4]{Sernesi}, the isomorphism classes of the deformations of the triple
$(X,\sL=\Oh_X(Z),\sigma)$ are the same as the isomorphism classes of the deformations of the pair $(X,Z=\{\sigma=0\})$.
Then, it is enough to apply Theorem \ref{main theorem of triple}.
 \end{proof}

 \begin{Remark}
 It is well known that the sheaf of DG-Lie algebras $ \Theta_X(-\log Z)$ controls the locally trivial deformations of the pair $(X,Z)$. In particular, it controls all infinitesimal deformations whenever $Z$ is smooth. The previous corollary generalises this result to the case of $Z$ any divisor in a smooth variety $X$ \cite[Section 8.1]{LMDT}.
 \end{Remark}

 \begin{Remark}
 If $\sL$ is a line bundle on a smooth Calabi--Yau variety $X$, then it is known that the infinitesimal deformations of the pair $(X,\sL)$ are unobstructed (see, for example, \cite{DM-linebundle,LP}).
Then, applying Corollary~\ref{cor H^1=0 allora forgetful smooth}, we recover the fact that if $H^1(X,\sL)=0$ then the infinitesimal deformations of the triples $(X,\sL, \sigma)$
are also unobstructed. This is equivalent to the deformations of the pair $(X,Z)$ being unobstructed and generalises \cite[Corollary 5.8]{donacoppie}, where the same result is proved only for smooth divisor (see also \cite{KKP}).
 \end{Remark}

\subsection*{Acknowledgements}
M.M.~is partially supported by the PRIN 20228JRCYB ``Moduli spaces and special varieties'', of ``Piano Nazionale di Ripresa e Resilienza, Next Generation EU''. D.I.~wishes to thank the Mathematical Department Guido Castelnuovo of Sapienza Universit\`a di Roma for the hospitality.
Both authors thank the anonymous referees for pointing out some minor mistake and for useful comments and suggestions on the paper.

\pdfbookmark[1]{References}{ref}
\LastPageEnding

\end{document}